\documentclass[11pt,reqno]{article}
\usepackage{amsmath,amsthm,amsfonts,amssymb,amscd}
\usepackage[latin1]{inputenc}
\usepackage{graphicx}

\def\leaderfill{\leaders\hbox to .8em{\hss .\hss}\hfill}
\def\_#1{{\lower 0.7ex\hbox{}}_{#1}}

\newtheorem{Theorem}{Theorem}

\newtheorem{Proposition}{Proposition}
\newtheorem{Lemma}{Lemma}

\theoremstyle{Definition}
\newtheorem{Definition}{Definition}

\theoremstyle{Remark}

\def\C{{\mathcal{C}}}

\def\G{{\mathcal{G}}}
\def\fa{{\mathcal{F}}}

\def\po{{\partial}}

\def\vr{{\varphi}}

\def\ov{\overline}

\def\re{{\mathbb{R}}}

\def\Diff{\operatorname{{Diff}}}
\def\sing{\operatorname{{sing}}}

\title{Foliations with Morse singularities}

\author{C\'esar Camacho   and Bruno Sc\'ardua}

\date{}

\begin{document}

\maketitle

\section{Introduction and main results}

The interplay between the topology of a closed manifold and the
combinatorics of the critical points of a real valued function of
class $C^2$ defined on the manifold is a well known fact of Morse
Theory. It is natural to expect a similar relationship for
foliated manifolds. This became evident for the first time with
the following result of G. Reeb ([10]), a consequence of his
Stability Theorem ([1],[9],[11]):
\begin{Theorem}
\label{Theorem : Reeb}  Let $M$ be a closed oriented and connected
manifold of dimension $m \ge 2$. Assume that $M$ admits a $C^1$
transversely oriented codimension one foliation $\fa$ with a non
empty set of singularities all of them centers. Then the singular
set of $\fa$ consists of two points and $M$ is homeomorphic to the
$m$-sphere.
\end{Theorem}

Later on Eells and Kuiper classified the closed manifolds
admitting a $C^3$ function with exactly three non-degenerated
singular points ([3]):
\begin{Theorem} [\cite{Ee-Ku1}]
\label{Theorem:eellskuiper} Let $M$ be a connected closed manifold
(not necessarily orientable) of dimension $m$. Suppose $M$ admits
a Morse function $f\colon M \to \re$ of class $C^3$ with exactly
three singular points. Then:
\begin{itemize}
\item[\rm(i)] $m \in \{2,4,8,16\}$ \item[\rm(ii)] $M$ is
topologically a compactification of $\re^m$ by an $\dfrac
m2$-sphere \item[\rm(iii)] If $m=2$ then $M$ is diffeomorphic to
$\re P(2)$. For $m \ge 4$\, $M$ is simply-connected and has the
integral cohomology structure of the complex projective plane
$(m=4)$, of the quaternionic projective plane $(m=8)$ on of the
Cayley projective plane $(m=16)$.
\end{itemize}
\end{Theorem}
We will call these manifolds $\it{Eells-Kuiper}$ manifolds.

In both situations we have a closed manifold endowed with a
foliation with Morse singularities where the number of centers is
greater than the number of saddles. In \cite{C-Sc} we proved that,
in the case that the manifold is orientable of dimension three,
this implies it is homeomorphic to the 3-sphere. The aim of this
paper is to consider the n-dimensional case. We proceed to define
the main notions we use. A codimension one foliation with isolated
singularities on a compact manifold $M$ is a pair $\fa = (\fa_0,
\sing\fa)$ where $\sing\fa \subset M$ is a discrete subset and
$\fa_0$ is a regular foliation of codimension one on the open
manifold $M\backslash\sing\fa$. We say that $\fa$ is {\it of
class} $C^k$ if $\fa_0$ is of class $C^k$, $\sing\fa$ is called
the singular set of $\fa$ and the {\it leaves} of $\fa$ are the
leaves of $\fa_0$ on $M\backslash\sing\fa$. A point $p \in
\sing\fa$ is a {\it Morse type} singularity if there is a function
$f_p\colon U_p \subset M \to \re$ of class $C^2$ in a neighborhood
of $p$ such that $\sing\fa \cap U_p = \{p\}$, $f_p$ has a
non-degenerate critical point at $p$ and the levels of $f_p$ are
contained in leaves of $\fa$. By the classical Morse Lemma there
are local coordinates $(y_1,\dots,y_m)$ in a neighborhood $U_p$ of
$p$ such that $y_j(p)=0$, $\forall\,j \in \{1,\dots,m\}$ and
$f(y_1,\dots,y_m) = f(p) - (y_1^2 +\cdots+ y_{r(p)}^2) +
y_{r(p)+1}^2 +\cdots+ y_m^2$\,. The number $r(p)$ is called the
Morse $\it{index}$ of $p$. The singularity $p$ is a {\it center}
if $r(p) \in \{0,m\}$ and it is a {\it saddle} otherwise. The
leaves of $\fa$ in a neighborhood of a center are diffeomorphic to
the $(m-1)$-sphere. Given a saddle singular point $p \in \sing\fa$
we have leaves of $\fa\big|_{U_p}$ that accumulate on $p$, they
are contained in the cone $\tau_{p}$: $y_1^2 +\cdots+ y_{r(p)}^2 =
y_{r(p)+1}^2 +\cdots+y_ {m}^2\neq 0$ and there are two
possibilities: either $r(p)=1$ or $m-1$ and then $\tau_{p}$ is the
union of {\it two} leaves of $\fa\big|_{U_p}$ , or $r(p)\neq 1$
and $ m-1$ and $\tau_{p}$ is a leaf of $\fa\big|_{U_p}$. Any leaf
of $\fa\big|_{U_p}$ contained in $\tau_{p}$ is called a {\it local
separatrix of $\fa$} at $p$, or a {\it cone leaf} at $p$ . Any
leaf of $\fa$ such that its restriction to $U_{p}$ contains a
local separatrix of $\fa$ at $p$ is called a {\it separatrix of
$\fa$ at $p$ }. A {\it saddle connection} for $\fa$ is a leaf
which contains local separatrices of two {\it different} saddle
points. A {\it saddle self-connection} for $\fa$ at $p$ is a leaf
which contains two different local separatrices of $\fa$ at $p$. A
foliation $\fa$ with Morse singularities {\it is transversely
orientable} if there exists a vector field $X$ on $M$, possibly
with singularities at $\sing\fa$, such that $X$ is transverse to
$\fa$ outside $\sing\fa$.

\begin{Definition}
\label{Definition:morsefoliation} {\rm A {\it Morse foliation}
$\fa$ on a manifold $M$ is a transversely oriented codimension one
foliation of class $C^2$ with singularities such that: (i) each
singularity of $\fa$ is of Morse type and (ii) there are no saddle
connections.} \end{Definition}

\noindent Basic examples of Morse foliations are given by the
levels of Morse functions $f\colon M \to \re$ of class $C^2$.
Therefore any manifold of class $C^2$ supports a Morse foliation,
i.e., the existence of a Morse foliation imposes no restriction on
the topology of the manifold. Nevertheless, there are restrictions
which come from the nature of the singularities of a Morse
foliation $\fa$ on $M$. Indeed, our purpose in this paper is to
show the following theorem:
\begin{Theorem}
\label{Theorem:main} Let $M$ be a compact connected manifold and
$\fa$ a Morse foliation on $M$ such that the number $k$ of
centers and the number $\ell$ of saddles in $\sing\fa$ satisfy $k
\ge \ell+1$. Then we have two possibilities:

{\rm (i)}\,$k=\ell+2$ and $M$ is homeomorphic to the $m$-sphere.

{\rm (ii)}\,$k=\ell+1$ and $M$ is an Eells-Kuiper manifold.
\end{Theorem}

\noindent{\bf Acknowledgement.} We want to thank Jos\'e Seade for
several valuable suggestions and the ICTP for the hospitality
during the elaboration of this work.

\section{Preliminaries}

Let us first fix the notation we use. Let $\fa$ be a Morse
foliation on a manifold $M$ of dimension $m \ge 3$. Given a center
singularity $p \in \sing\fa$ the nearby leaves of $\fa$ are
compact diffeomorphic to $S^{m-1}$. Since $m-1 \ge 2$ any such
leaf $L$ has trivial holonomy group and therefore by the Local
 Stability theorem of Reeb (\cite{C-LN}, \cite{Godbillon}) there is a fundamental
system of open neighborhoods $V$ of $L$ such that the restriction
$\fa|_V$ is equivalent to a product foliation $\G$ on $L \times
(-1,1)$ whose leaves are of the form $L \times \{t\}$, $t \in
(-1,1)$. We introduce therefore the {\sl open\/} subset $\C(\fa)$
as the union of centers in $\sing\fa$ and leaves of $\fa$
diffeomorphic to $S^{m-1}$. The {\it basin} of a center $p\in\sing
\fa$ is the connected component  $\C_p(\fa)$  of $\C(\fa)$ that
contains  $p$. We have the following basic lemma:

\begin{Lemma}
\label{Lemma:basic} Given centers $p,q \in \sing\fa$ the sets
$\C_p(\fa)$ and $\C_q(\fa)$ are open in $M$ and $\C_p(\fa) \cap
\C_q(\fa) \ne \emptyset$ if and only if $\C_p(\fa) = \C_q(\fa)$.
Moreover we have $\C_p(\fa) = M$ if and only if $\po\C_p(\fa) =
\emptyset$ and in this case $M$ is homeomorphic to $S^m$ provided
that $M$ is orientable.
\end{Lemma}

\noindent In particular  $M$ is homeomorphic to $S^{m}$ or
$\po\C_p(\fa)$ contains some saddle singularity.

 In order to fix notations we shall now introduce the notion of
holonomy group
of an invariant subset of codimension one. We will consider two
notions of holonomy. When we refer to the {\it holonomy of a leaf}
$L$ of $\fa$ we mean the holonomy group of $L$ as a leaf of
$\fa_0$ on $M\setminus \sing\fa$. On the other hand, the notion of
holonomy can be extended to invariant subsets of the form
 $S=\tau\cup \{p\}$, \, $ p \in \sing\fa$ and
$\tau$ is either a cone leaf or a union of two cone leaves. Notice
that, in a small neighborhood of $p$, $\tau$ can consist of two
components $\tau_1$ and $\tau_2$, and that this can only happen if
$r(p)= 1$ or $m-1$. In this case $S$ locally divides the manifold
into three connected components. One of them, say $R_3$, is the
union of (regular) leaves which are hyperboloids of one sheet, and
the others, say $R_1$ and $R_2$, are the union of one of the
connected components of hyperboloids of two sheets (we can think
of $R_1$ as the region surrounded by $\tau_1$ and $R_2$ is the
region surrounded by $\tau_2$). Let $\gamma\colon [0,1] \to S$ be
a path on $S$  which passes through the singularity $p$ (from
$\tau_1$ to $\tau_2$). In this case the holonomy along $\gamma$
can be defined in the usual manner (lifting paths to leaves) on
$R_3$, however, there is no canonical extension of this holonomy
to the other side in general. Thus we adopt the following notion
of holonomy. Fix a neighborhood $U$ of $p\in \sing\fa$ where $\fa$
is given by a Morse function $f$ with a single singularity at $p$.
Let $\gamma\colon [0,1] \to S$ be a piecewise smooth    path (as a
map $\gamma\colon [0,1]\to M$). Let $T_0$ and $T_1$ be local
transversals to $\fa$ at $\gamma(0)$ and $\gamma(1)$ respectively.
The {\it holonomy} along $\gamma$ will be the mapping which
assigns $t\in T_0$ to $f^{-1}(f(t))\cap T_1 \in T_1$. This
holonomy map is well-defined even if $\gamma$ is not contained in
$\{p\}\cup \tau_1$. \vglue.1in

 Next we study the possible
intersections for the boundaries of two basins.

\begin{Lemma}
\label{Lemma:topologyseparatrix} Suppose $p_1,p_2 \in \sing\fa$
are distinct centers such that $\po\C_{p_1}(\fa) \cap \po\C_{p_2}(\fa)
\ne \emptyset$. Then we have the following mutually exclusive possibilities:
\begin{itemize}
\item[\rm(i)] $\po\C_{p_1}(\fa) = \po\C_{p_2}(\fa)$ and so $M =
\ov{\C_{p_1}(\fa)} \cup \ov{\C_{p_2}(\fa)}$.
\item[\rm(ii)]
 $\po\C_{p_1}(\fa) \ne \po\C_{p_2}(\fa)$ and there is a saddle point $q
\in \po\C_{p_1}(\fa) \cap \po\C_{p_2}(\fa)$ with Morse index $1$ or
$m-1$, without self connection.
\end{itemize}
\end{Lemma}

\begin{proof} Since $\po\C_{p_1}(\fa) \cap
\po\C_{p_2}(\fa) \ne \emptyset$, by Lemma 1 there is a saddle singular
point $q \in\po\C_{p_1}(\fa) \cap \po\C_{p_2}(\fa)$ .
 If the Morse index of $q$ is different from $1$
and $m-1$ then in  suitable local coordinates $(y_1,\dots,y_m)$
we have $q = (0,\dots,0)$ and the local separatrix $\tau_{q}$ through
$q$ is given by $y_1^2 +\cdots+ y_r^2 = y_{r+1}^2+\cdots+
y_m^2 \neq 0$ where $r \notin \{1, m-1\}$. In particular $ \tau_{q} $
is connected.

Thus, if $C$ is the separatrix of $\fa$ at $q$ we have $\po\C_{p_1}(\fa) =
\po\C_{p_2}(\fa) = \ov C = C \cup \{q\}$ and we are in case (i).

\noindent In case the index of $q$ is 1 or $m-1$, and $C$ is a
self-connection at $q$, then $\po\C_{p_1}(\fa) = \po\C_{p_2}(\fa) =
\ov C$ and we are again in case (i). The remaining case is index
of $q$ is $1$ or $m-1$ and $q$ has
no self-saddle connection. Consider local coordinates
$(y_1,\cdots,y_m)$ with $q=(0,\cdots,0)$ and $\fa$ given
by the levels of the function $-y_1^2+y_2^2+\cdots+y_m^2$.
The level zero of this function bounds the regions $R_1$:
$y_1<0$, $y_1^2 > y_2^2+\cdots+y_m^2 $, $R_2$: $y_1>0$, $y_1^2>
y_2^2+\cdots+y_m^2$ and $R_3$ : $y_1^2 < y_2^2+\cdots+
y_m^2$. Then $\C_{p_i}(\fa)\cap R_3=\emptyset$, i=1,2, because
otherwise we would have a saddle self-connection at $q$. On the
other hand $\C_{p_1}(\fa) \cap R_1 \neq\emptyset$ implies
$\C_{p_1}(\fa)\cap R_2=\emptyset$ by the same reason. Therefore
$\po\C_{p_1}(\fa)\neq\po\C_{p_2}(\fa)$. This proves (ii).\end{proof}

\begin{Proposition}
\label{Proposition:Miseellskuiper} Let $\fa$ be a Morse foliation
on a closed connected manifold $M$ of dimension $m \ge 3$. Assume
that $k=2$ and $\ell=1$, i.e., $\fa$ has exactly two centers and
one saddle singularity. Then $M$ is homeomorphic to an
Eells-Kuiper manifold.
\end{Proposition}

\begin{proof} We shall first prove that the
nonsingular foliation $\fa_0 = \fa\big|_{M_0}$ on
$M_0 = M\backslash\sing\fa$ is a proper stable foliation.
 There are several equivalent conditions that
define a stable foliation (\cite{Godbillon}). We shall prove that
given any leaf $L_0$ of $\fa_0$ there is a fundamental system of
open neighborhoods  of $L_0$ in $M_0$ saturated by $\fa_{\!0}$. By
hypothesis we have $M = \ov{\C_{p_1}(\fa)} \cup \ov{\C_{p_2}(\fa)} =
\C_{p_1}(\fa) \cup \C_{p_2}(\fa) \cup C\cup \{q\} $ where $\po\C_{p_1}(\fa)$
=$\po\C_{p_2}(\fa)$=$\ov{C}$. Thus if $r(q)=1$ or $m-1$ then $C$
is a self connection.We proceed to show that this cannot occur.
Suppose that $\C_{p_1}(\fa)\cap R_3\neq \emptyset$, then $\C_{p_2}(\fa)$
will have nonempty intersection with $R_1$ and with $R_2$.
Taking a small closed ball $\ov{B_q(t)}$ of radius $t > 0$ centered at
$q$ then for leaves
$L_i $ of $ \fa$, $L_i \subset \C_{p_i}(\fa)$ for $i \in \{1,2\}$,
close enough to $C$,  we have that
the intersection $L_1 \cap (M\backslash\ov{B}_q(t))$ is a union of two
disjoint $(m-1)$-discs. Moreover $L_2 \cap (M\backslash\ov{B}_q(t))$
is the complement of two disjoint $(m-1)$-discs in an $(m-1)$-
sphere. Since both manifolds $L_1 \cap (M\backslash\ov{B}_q(t))$ and
$L_2 \cap (M\backslash\ov{B}_q(t))$ are homeomorphic to
$C\backslash (C\cap\ov{B}_q(t))$ then we obtain a contradiction.
Thus, $r\neq 1,m-1$ and
$C\cap\ov{B}_q(t)$ is connected for t small.
\noindent Given a leaf $L_0$ of $ \fa$ we have two
possibilities, either $L_0 \subset\C_{p_i}(\fa)$ for some $i \in
\{1,2\}$ and $L_0$ is homeomorphic to $S^{m-1}$, or $L_0 =C$.
In case $L_0$ is in $\C_{p_i}(\fa)$ then, by
the Reeb stability theorem, $L_0$ has a fundamental system of
saturated neighborhoods consisting of compact leaves. This shows
that the leaves of $\fa_0$ in $M_0\backslash C$ are stable. It
suffices to show that $C$ is also a stable leaf.

\noindent We claim that the holonomy group of $C\cup \{q\}$ is a
finite group conjugated to a subgroup of
$\Diff(\re,0)$. Indeed, $C\backslash\ov{B_q(t)}$ is a disc and
therefore simply-connected; on the other hand in $\ov{B_q(t)}$ the
foliation has a first integral as $f = -( y_1^2 +\cdots+ y_r^2) +
y_{r+1}^2 +\cdots+ y_m^2$ so that the holonomy group of $\fa
\cap \ov{B_q(t)}$ is finite.Since $C\cup{\{q\}}$ has a finite holonomy
group, which is a subgroup of $\Diff(\re,0)$, the holonomy group of
$C\cup\{q\}$ is either trivial or, in case $\fa$ is not transversely
orientable, it has order 2. By the classical argument on stability
of leaves we conclude that finite holonomy implies that the leaf $C$
is stable and $\fa_0$ is stable in $M_0$.

 \noindent Since
$\fa_0$ is stable in $M_0$\,, the leaf space $M_0\big/\fa_0 =:
\frak{X}_{\fa_0}$ is Hausdorff and therefore it is a 1-manifold.
The choice of a differentiable submersion $\frak{X}_{\fa_0} \to
\re$ gives then a differentiable first integral $F_0\colon M_0 \to
\re$ for $\fa_0$\,. Clearly $F_0$ can be modified in order to
admit a differentiable (radial) extension to the center
singularities $p_1,p_2 \in \sing\fa$. It remains to show that
$F_0$ can be modified in order to admit a differentiable extension
to $q$. This is a consequence of the fact that by the triviality
of the holonomy group of $C\cup{q}$ we can extend the local first
integral $f = - \sum\limits_{j=1}^r y_j^2 + \sum\limits_{k=r+1}^m
y_k^2$ from a ball $B_q(t)$ to a neighborhood $T$ of $C\cup\{q\}$ in
$M$ in such a way that $\po T$ is a union of leaves of $\fa$, each
leaf diffeomorphic to an $(m-1)$-sphere and contained in some basin
$\C_{p_i}(\fa)$.

\noindent Thus we have proved that $\fa$ is given by a Morse
function $F\colon M \to \re$ and therefore, by Eells-Kuiper
Theorem~\ref{Theorem:eellskuiper}, $M$ is an Eells-Kuiper
manifold.
\end{proof}

For any $r>0$ we will write $B(r)=D(r)\times I(r)\subset \mathbb
R^{m-1} \times \mathbb R$, where $D(r)$ and $I(r)$ are closed
discs of radius $r$ centered at zero. The foliation on $B(r)$
given by the submersion $(x,t)\mapsto t$ will be denoted by
$\mathcal H$. Let $p\in M$  be a center and let $q \in \partial
\mathcal C_p(\fa)$ be a saddle point of $\fa$. We will say that
$p, q$ {\it form a trivial pairing} if there are open
neighborhoods $V\supset V^\prime \supset \ov{\mathcal
C_p(\fa)}, p, q \in V^\prime$ and a diffeomorphism $\vr \colon \ov
V \to B(1)$, onto $B(1)$, such that $\vr (\ov V^\prime)=B(1/2)$
and $\fa\big|_{\ov V\setminus V^\prime}=\vr^* \mathcal H$.

\begin{Lemma}
Suppose $p_1, p_2 \in M$ are two different centers such that
$\partial \mathcal C_{p_1}(\fa)\cap \partial \mathcal
C_{p_2}(\fa)\ne \emptyset$ and let $q$ be the saddle point
contained in this intersection. Assume that the index of $q$ is
one and that there is no saddle self-connection at $q$. Then,
either $p_1, q$ or $p_2, q$ form a trivial pairing.
\end{Lemma}

\begin{proof}
In a neighborhood $q\in U$ there are local coordinates
$(y_1,...,y_m)$ such that $q=(0,...,0)$ and the leaves of
$\fa\big|_U$ are given by the levels of the function
$f(y_1,...,y_m)=-y_1 ^2 +(y_2 ^2 +...+y_m^2)$ on $\mathbb R^m$.
As before the cone $-y_1 ^2 +(y_2 ^2 +...+y_m^2)=0$ divides $U$ in three
regions. The regions $R_1$ and $R_2$ are defined by $R_1$:
$y_1<0$, $y_1 ^2 >y_2
^2 +...+y_m^2$ , $R_2$: $y_1>0$,  $y_1 ^2 >y_2
^2 +...+y_m^2$ and the region $R_3$ by $y_1 ^2 <y_2 ^2
+...+y_m^2$. The cone leaves $\tau_1$ and $\tau_2$ of $\fa\big|_U$ bound
$R_1$ and $R_2$ respectively and since there is no self-connection
at $q$ there are different leaves of $\fa$, $\ell_1$ and $\ell_2$,
such that $\ell_1 \supset \tau_1$ and $\ell_2 \supset \tau_2$. Since we
have two center basins $\mathcal C_{p_1}(\fa)$ and $\mathcal
C_{p_2}(\fa)$ and three regions $R_1, R_2, R_3$ then some
$\mathcal C_{p_i}(\fa)$ will intersect $R_1$ or $R_2$. Suppose
that $\mathcal C_{p_1}(\fa)\cap R_1 \ne \emptyset$. Then $\mathcal
C_{p_1}(\fa)\cap R_2 = \emptyset$ and $\mathcal C_{p_1}(\fa)\cap
R_3=\emptyset$ because both $R_2$ and $R_3$ have $\ell_2$ in their
boundary and if either $\mathcal C_{p_1}(\fa)\cap R_2\ne
\emptyset$ or $\mathcal C_{p_2}(\fa)\cap R_3\ne \emptyset$ this
would imply a saddle self-connection at $q$. Thus $\partial
\mathcal C_{p_1}(\fa)= \ov \ell_1$. We claim that $\ov \ell_1=
\ell_1 \cup \{q\}$. Indeed, if on the contrary, there is a regular
point of $\fa$, $s\in \ov \ell_1 \setminus \ell_1$, and $S_s$
denotes an arbitrarily small cross section to $\fa$ centered at
$s$, then the number of points of intersection of $\ell_1$ with
$S_s$, $n(\ell_1, S_s)$, is infinite. On the other hand, by Reeb's
theorem, given any local transverse section $S$ to $\fa$, the
number of points of intersection $n(\ell, S)$, of a leaf $\ell
\subset \mathcal C_{p_1}(\fa)$ with $S$, is locally constant.
Since $\ell_1$ is approached by leaves in $\mathcal C_{p_1}(\fa)$
we would obtain leaves $\ell(k)\subset \mathcal C_{p_1}(\fa)$ with
$n(\ell(k), S_s)\to \infty$ as $k\to \infty$, which is a
contradiction. Therefore $\partial \mathcal C_{p_1}(\fa)= \ell_1
\cup \{q\}$.

We take $U$ small enough so that $\ov \ell_1 \cap U=\ov \tau_1$.
Thus, for any leaf $L\subset \mathcal C_{p_1}(\fa)$ close enough
to $\ov \ell_1$, $L\cap U$ is connected and $f\big|_{L\cap
U}=\delta <0$, a constant. Write $L=L_\delta$. Then, as $\delta
\to 0$ $L_\delta\setminus U$ approaches $\ov \ell_1 \setminus \tau_1$
and this implies that $\ov \ell_1 \setminus \tau_1$ is homeomorphic
to $L_\delta \setminus U$, i.e., to an $(m-1)-$disc. Therefore
$\ov \ell_1$ is homeomorphic to $S^{m-1}$. Moreover for $\epsilon
>0$ small enough each leaf $f^{-1}(\delta)$ of $\fa\big|_{R_1\cup
R_3}$, $-2\epsilon \leq \delta \leq 2\epsilon$, bounds an
$(m-1)-$disc $D_\delta$ close to $\ell_1\setminus \tau_1$, with
$D_0=\ell_1 \setminus \tau_1$. The union $T_{2\epsilon}=
\bigcup\limits_{-2\epsilon \leq \delta \leq 2 \epsilon} \ov
D_\delta$ is a trivially foliated neighborhood of $\ell_1
\setminus \tau_1$. We can extend the function $f$ to $U \cup
T_{2\epsilon}$ by writing $f \big|_{\ov D _\delta}=\delta$. We
define a saturated neighborhood $V_0$ of $\ov \ell_1 \cup \tau_2$ by
$V_0=f^{-1}([-\epsilon, \epsilon]) \cup \mathcal C_{p_1}(\fa)$ and
define $g=f\big|_{{V_0}\setminus\C_{p_1}(\fa)}$.
Consider now a Riemannian metric defined on $M$ and a normal
vector field to $\fa$ that in $U$ takes the form $\mathcal N=-
y_1\frac{\partial}{\partial y_1} +y_2 \frac{\partial}{\partial
y_2} +...+ y_m\frac{\partial}{\partial y_m}$.
For $a>0$ small enough the submanifold $e=(y_1=a)\cap \tau_2$ is
well defined and diffeomorphic to $S^{m-2}$. Let $\Sigma$ be
a cross section to $\fa$, over e, i.e. $\Sigma\cap \tau_2=e$,
$\Sigma$ contained in $V_0$ and invariant by $\mathcal N$. We
take $\Sigma$ diffeomorphic to $e\times [-\epsilon, \epsilon]$
by means of a map that takes each $e\times{\delta}$, $-\epsilon
\leq \delta\leq \epsilon$, to the leaf $(f=\delta)\cap\Sigma$
of $\fa\big|_\Sigma$.

Consider the region $V\subset V_0$ bounded by $(g=-\epsilon)\cup
(g=\epsilon)$ and $\Sigma$ and define a neighborhood
$\partial V \subset W \subset V$, as $W=(-\epsilon\leq
g \leq -\epsilon/2) \cup ( \epsilon/2 \leq g \leq
\epsilon)\cup N$ where $N$ is a neighborhood of
$\Sigma\subset V$ invariant by $\mathcal N$.

The foliation $\fa\big|_W$ is trivial in the sense that on
$(-\epsilon\leq g\leq -\epsilon/2)\cup (\epsilon/2\leq g\leq\epsilon)$
the leaves of $\fa$ are levels of $g$, diffeomorphic to $D(1)$
and on $N$ the leaves of $\fa$ are levels of $g$
$(g=\delta)$, $-\epsilon/2\leq\delta\leq\epsilon/2$
diffeomorphic to $D(1)\setminus D(1/2)$.
Moreover in $W$ the foliation $\fa$ and the trajectories of $\mathcal
N$ are everywhere transverse. Thus a diffeomorphism
$\vr \colon W \to
B(1)\setminus B(1/2)$ can be easily constructed by sending leaves
of $\fa\big|_W$ at the level $(g=\alpha \epsilon)$ onto leaves of
$\mathcal H\big|_{B(1)\setminus B(1/2)}$ at the level $(t=\alpha)$
and orbits of $\mathcal N\big|_W$ to orbits of
$\frac{\partial}{\partial t}$. Then $\vr$ is extended to $V$.
\end{proof}

\section{Proof of the Theorem}

Now we prove Theorem~\ref{Theorem:main}.  We will proceed by
induction on the number $\ell$ of saddle singularities. If
$\ell=0$ then Reeb's Theorem applies and $M$ is homeomorphic to
$S^m$. Assume now that $\ell \ge 1$ and the result has been proven
for foliations with at most $\ell-1$ singularities of saddle type.
We have $k \ge \ell + 1$. Thus $k \ge 2$. Suppose that $M$ is not
homeomorphic to $S^m$. Then for each center $p \in \sing\fa$ there
must be a saddle $q(p) \in \po\C_p(\fa)$ .Since $k \ge \ell + 1$
and $k \ge 2$ there are two centers $p_1$, $p_2$ such that $q(p_1)
= q(p_2)$, i.e., there is a saddle $q$ such that $q \in
\po\C_{p_1}(\fa) \cap \po\C_{p_2}(\fa) $ and by Lemma 3 we have
either $M =\ov{\C_{p_1}(\fa)} \cup \ov{\C_{p_2}(\fa)}$ or $q$ has
index 1 or $m-1$ and is not self-connected.

\noindent In case $M = \ov{\C_{p_1}(\fa)} \cup \ov{\C_{p_2}(\fa)}$
then clearly $\C_{p_i}(\fa) \cap \sing\fa = \{p_i\}$ , i=1,2.

\noindent Thus $\sing\fa = \{p_1, p_2, q\}$ and by the
Proposition $M$ is an
Eells-Kuiper manifold. In case $q$ has index 1 or $m-1$ and no
self-connection then by Lemma 3 we can
eliminate one saddle and one center replacing $\fa$ by a Morse
foliation $\fa_1$ on $M$ with a number $k_1$of centers and $\ell_1$
of saddles given by $k_1 = k-1$ and $\ell_1 =\ell-1$. Therefore $k_1
\ge \ell_1+1$ and $\ell > \ell_1 \ge 0$. By the induction
hypothesis $M$ is homeomorphic to $S^m$ or to an Eells-Kuiper
manifold. This proves the theorem.\qed

\bibliographystyle{amsalpha}

\begin{thebibliography}{10}
\frenchspacing




\bibitem{C-LN} C. Camacho, A. Lins Neto:
Geometric  Theory of Foliations, Birkh\"auser; 1985.

\bibitem{C-Sc} C. Camacho, B. Sc\'ardua: {\it On
codimension one foliations with Morse singularities on
three-manifolds}, to appear in Topology and its Applications.


\bibitem{Ee-Ku1}   J. Eells, N. Kuiper: {\it
Manifolds which are like projective planes}; Pub. I.H.E.S., t.14
(1962), p. 5-46.

\bibitem{Ee-Ku2} J. Eells, N. Kuiper: {\it Closed manifolds
which admit nondegenerate functions with three critical points};
Indag. Math. 23 (1961), p. 411-417.


\bibitem{Godbillon} C. Godbillon: Feuilletages.  \'Etudes g\'eom\'etriques.
With a preface by G. Reeb. Progress in Mathematics, 98.
Birkh\"auser Verlag, Basel, 1991.


\bibitem{H1}  A. Haefliger;
{\em Structures feuillet\'{e}es et cohomologie \`{a} valeur dans
un faisceau de group\"oides.}  Comment. Math. Helv. 32 1958
248--329.

\bibitem{H2}  A. Haefliger; {\it Vari\'et\'es
feuillet\'ees}, An. Scuola Norm. Sup. Pisa, 3, vol. 16, 1962,
367-397.




\bibitem{M} J. Milnor: Morse Theory,
Annals of Mathematics Studies Number 51, Princeton, N.J., 1963.


\bibitem{R1}  G. Reeb; {\it Sur certaines propri\'et\'es
topologiques des vari\'et\'es feuillet\'ees}; Actualit\'es Sci.
Ind., Hermann, Paris, 1952.

\bibitem{R2}  G. Reeb; {\it Sur les points singuliers
d'une forme de Pfaff compl\'etement int\'egrable ou d'une fonction
num\'erique}. C.R.A.S. Paris 222, 1946, p.847-849.


\bibitem{Wa} E. Wagneur;
{\it Formes de Pfaff à singularités non dégénérées}. Annales de
l'institut Fourier, 28 no. 3 (1978), p. 165-176


\end{thebibliography}

\begin{tabular}{ll}
C\'esar Camacho  & \qquad  Bruno Sc\'ardua\\
IMPA-Estrada D. Castorina, 110 & \qquad Instituto de Matem\'atica\\
Jardim Bot\^anico  & \qquad Universidade Federal do Rio de Janeiro\\
Rio de Janeiro - RJ  & \qquad  Caixa Postal 68530\\
CEP. 22460-320   & \qquad 21.945-970 Rio de Janeiro-RJ\\
BRAZIL &  \qquad BRAZIL
\end{tabular}

\end{document}